\documentclass[12pt]{article}
\usepackage{latexsym}
\newtheorem{theorem}{Theorem}
\newtheorem{lemma}{Lemma}
\def\X{X}
\def\E{E}
\def\O{{\rm O}}
\def\p{p\,}
\def\mod{\mathop{\rm mod}}
\def\Re{\mathop{\rm Re}}

\begin{document}
\title{Regularly Spaced Subsums of Integer Partitions}
\author{E. Rodney Canfield \\ \small{Department of Computer
Science, University of Georgia}\\ \small{Athens, GA 30602}\\\texttt{\small<erc@cs.uga.edu>}
\and
Carla D. Savage
\thanks{Research supported by NSA grant MDA 904-01-0-0083} \\ \small{Department of Computer
Science, North Carolina State University}\\\small{ Raleigh, NC 27695-8206}\\\texttt{\small<savage@unity.ncsu.edu>}
\and
Herbert S. Wilf \\\small{Department of Mathematics, University of
Pennsylvania}\\\small{Philadelphia, PA 19104-6395}\\\texttt{\small<wilf@math.upenn.edu>}
\date{}}
\maketitle

\begin{abstract}
For integer partitions $\lambda :n=a_1+...+a_k$, where $a_1\ge a_2\ge \dots\ge a_k\ge 1$, we study the sum $a_1+a_3+\dots$ of the parts of odd index. We show that the average of this sum, over all partitions $\lambda$ of $n$, is of the form $n/2+(\sqrt{6}/(8\pi))\sqrt{n}\log{n}+c_{2,1}\sqrt{n}+O(\log{n}).$ More generally, we study the sum $a_i+a_{m+i}+a_{2m+i}+\dots$ of the parts whose indices lie in a given arithmetic progression and we show that the average of this sum, over all partitions of $n$, is of the form $n/m+b_{m,i}\sqrt{n}\log{n}+c_{m,i}\sqrt{n}+O(\log{n})$, with explicitly given constants $b_{m,i},c_{m,i}$. Interestingly, for $m$ odd and $i=(m+1)/2$ we have $b_{m,i}=0$, so in this case the error term is of lower order. The methods used involve asymptotic formulas for the behavior of Lambert series and the Zeta function of Hurwitz.

We also show that if $f(n,j)$ is the number of partitions of $n$ the sum of whose parts of even index is $j$, then for every $n$, $f(n,j)$ agrees with a certain universal sequence, Sloane's sequence \texttt{\#A000712}, for $j\le n/3$ but not for any larger $j$. 
\end{abstract}

\newpage

\section{Introduction}
If
\[\lambda:\ n=a_1+a_2+\dots +a_k \qquad (a_1\ge a_2\ge\dots\ge a_k\ge 1)\]
is a partition of $n$, we will refer to $a_1,a_3,a_5,\dots$ (resp.
$a_2,a_4,a_6,\dots$) as the parts of \textit{odd index} (resp. of
\textit{even index}) of the partition $\lambda$, and we define the two associated
partitions
\[\lambda_o:\ a_1+a_3+a_5+\dots,\qquad \lambda_e:\
a_2+a_4+a_6+\dots. \]
Our interest in these partitions was stimulated by two developments. First,
Astrid Reifegerste \cite{re} has shown how the sign of a permutation
$\sigma$ can be deduced from the partitions $\lambda_e,\lambda_o$ that belong to the
shape of the tableaux that are induced by $\sigma$ under the RSK
correspondence.
Second, we have observed the following interesting phenomenon. Suppose we
define $f(n,j)$ to be the number of partitions of $n$ such that $|\lambda_e|=j$.
Then the sequence $\{f(20,j)\}$ begins as
\[ 1,2,5,10,20,36,65,109,167,170,42,0,0,0,\ldots \]
and the sequence $\{f(25,j)\}$ begins as
\[1,2,5,10,20,36,65,110,185,297,443,512,272,0,0,0,\ldots .\]
As $n$ increases, the values of $f$ seem to be approaching the values of a
certain universal sequence which begins as
\[1,2,5,10,20,36,65,110,185,300,481,752,1165,1770,2665,3956,\dots .\]
Reference to the database \cite{is} quickly reveals that sequence
\#\texttt{A000712} is identical to the above as far as the computations go.
That sequence is described as the number of ``colored partitions,'' meaning
that we color the parts of some partition in two colors, and regard two such
partitions as being
the same colored partition if for each part $i$, the number of copies
of $i$ of a given color is the same in both partitions.
Another way to describe this is to say that
the $(n+1)^{\mathrm{st}}$ member of sequence \#\texttt{A000712} counts
ordered pairs $(\alpha,\beta)$ of integer partitions such that
$|\alpha|+|\beta|=n$. This latter viewpoint shows that
\begin{equation}\label{eq:form}\#\mathtt{A000712}(n+1)=\sum_jp(j)p(n-j)
\qquad (n=0,1,2,3,\dots ),\end{equation}
where $p$ is the usual partition function.
Closer inspection of that sequence and of our computations suggests that the
following more precise statement might be true.
\begin{theorem}
\label{th:one}
For each $n=0,1,2,\dots$, the number of partitions of $n$ the sum of whose
parts of even index is $j$ is equal to the right side of (\ref{eq:form})
above, for $0\le j\le n/3$, and that bound for $j$ is best possible.
\end{theorem}
We prove this theorem in Section 2.
Next we consider the relative contributions of the parts of even and of odd
indices of a typical partition. Evidently $|\lambda_o| \ge |\lambda_e|$
always, so $|\lambda_e|\leq |\lambda|/2$. However one might expect that the
parts of even index, even though they contribute less than $|\lambda|/2$, might
not contribute very much less than that, so on the average the two halves
might contribute asymptotically the same amounts. 
And what about a modulus other than
2?  For modulus $m \geq 1$ and for
$1\leq i\leq m$,
let $\X_{m,i}(\lambda)$ be the sum of those
parts in partition $\lambda$ whose index $j$ is congruent to $i$
mod $m$:
\begin{equation}
\label{eq:Xdef}
\X_{m,i}(\lambda)=\sum_{j:j\equiv i (\mod m)}\lambda_j.
\end{equation}
One might expect that on the average, each $X_{m,i}$ is about $n/m$, and 
this is indeed the case.  More specifically,
let $C=\pi\sqrt{2/3}$; we will prove in Section 4:
\begin{theorem}
\label{th:three}
For fixed integers $m\ge 1$ and $i$, there exists a constant
$c_{m,i}$ such that
\begin{equation}
\label{eq:Th31}
\E(\X_{m,i}) - \frac{n}{m} = \frac{m+1-2i}{2Cm}\sqrt{n}\log{n}
+ c_{m,i}\sqrt{n} + \O(\log{n}).
\end{equation}
The constants $c_{m,i}$ are given by
\begin{equation}
\label{eq:Th32}
c_{m,i} = \frac{\bigl(\gamma+\log(2/C)\bigr)\bigl(m+1-2i\bigr)}{Cm} 
+ \frac{2}{Cm}\sum_{\ell=1}^{m-1}{
                  \omega^{-\ell(i-1)} \over 1-\omega^{\ell}
                                 }\log(1-\omega^{\ell}),
\end{equation}
where $\omega=e^{2\pi\sqrt{-1}/m}$, and  
$\gamma=0.5772\cdots$ is the well-known Euler constant.
\end{theorem}
An interesting thing to note about
this result is that when the modulus $m$ is odd and $i=(m+1)/2$,
then the term on the right of equation (\ref{eq:Th31}) having magnitude $\sqrt{n}\log{n}$
disappears, and the difference $\E(\X_{m,i})-n/m$ has
magnitude  $\sqrt{n}$.  For example, when $m=3$
and $i=2$, so we are taking the expected value of $\lambda_2
+\lambda_5+\cdots$, the theorem asserts
$$
\E(X_{3,2}) - \frac{n}{3} = -\frac{\sqrt{2}}{9}\sqrt{n} + \O(\log{n}).
$$

\section{Proof of Theorem \ref{th:one}}
Recall that we have defined $f(n,j)$ to be the number of partitions of $n$
the sum of whose parts of even index is $j$.
Let $p(n,j)$ be the number of partitions of $n$ with at most $j$ parts.
We first obtain an explicit formula for $f$, viz.,
\begin{equation}
\label{eq:fform}
f(n,j)=\sum_{i\le j}p(i)p(j-i,n-2j),
\end{equation}
by describing a bijection between the set of partitions $\lambda$ of $n$
with $|\lambda_e|=j$ and
 the set of pairs of partitions
$(\alpha,\beta)$ where $|\alpha|+|\beta|=j$ and $\beta$ has at most $n-2j$
parts.

Let $\lambda =(a_1,a_2, \ldots, a_r)$ be a partition of $n$ with
$j= |\lambda_e|$.  Map $\lambda$ to the pair of partitions $(\alpha,\beta)$
with
\[
\alpha = (1^{a_2-a_3}2^{a_4-a_5}3^{a_6-a_7} \cdots), \ \ \ \ \
\beta = (1^{a_3-a_4}2^{a_5-a_6}3^{a_7-a_8} \cdots),
\]
where we let $a_i=0$ if $i>r$.
Then
\begin{eqnarray*}
|\alpha|+|\beta| & =&
 \sum_{i\geq 1}i((a_{2i}-a_{2i+1})+(a_{2i+1}-a_{2i+2})) \\
& = &
 \sum_{i\geq 1}i(a_{2i}-a_{2i+2}) \\
& = &
 \sum_{i\geq 1}a_{2i} = |\lambda_e|= j, \\
\end{eqnarray*}
and the number of  parts of $\beta$ is
\[
(a_3-a_4) +(a_5-a_6) + \cdots \leq |\lambda_o|-|\lambda_e|=(n-j)-j=n-2j.
\]

To show this mapping is a bijection,
assume $n$ and $j$ are fixed and $(\alpha,\beta)$ is
a pair of partitions satisfying $|\alpha|+|\beta|=j$ and the number of
parts of $\beta$ is $b \leq n-2j$.  Let $\delta = n-2j-b$.
For a partition $\lambda$ and positive integer $i$, define $m_{\lambda}(i)$
to be the multiplicity of part $i$ in $\lambda$.  The inverse mapping
sends $(\alpha,\beta)$ to the partition
$\lambda = (a_1, a_2, \ldots )$ defined by
\[
a_1=a_2+ \delta; \ \ \ \ 
a_i = 
\sum_{k \geq \lceil i/2 \rceil} m_{\alpha}(k) +
\sum_{k \geq \lfloor i/2 \rfloor} m_{\beta}(k),  \ \ \ \ \ {i \geq 2}.
\]
Clearly $\lambda$ is a partition.
Note that
\[
a_1 = \sum_{k \geq 1} m_{\alpha}(k) +
\sum_{k \geq 1} m_{\beta}(k)
+n-2j
-\sum_{k \geq 1} m_{\beta}(k)= n-2j+\sum_{k \geq 1} m_{\alpha}(k).
\]
The weight of $\lambda_e$ is
\begin{eqnarray*}
|\lambda_e| &  = & \sum_{i \geq 1} a_{2i} =
\sum_{i \geq 1}\left(
\sum_{k \geq  i } m_{\alpha}(k) +
\sum_{k \geq  i } m_{\beta}(k)
\right)\\
& = &
\sum_{k \geq  i} km_{\alpha}(k) +
\sum_{k \geq  i} km_{\beta}(k)\\
& = & |\alpha|+|\beta|=j,
\end{eqnarray*}
and the weight of $\lambda$ is
\begin{eqnarray*}
|\lambda| & = & a_1 + \sum_{i \geq 1} \left(
\sum_{k \geq \lceil i/2 \rceil} m_{\alpha}(k) +
\sum_{k \geq \lfloor i/2 \rfloor} m_{\beta}(k)
\right)\\
& = & a_1 + (2|\alpha|- \sum_{k \geq 1} m_{\alpha}(k) )+2|\beta|\\
& = & n-2j+\sum_{k \geq 1} m_{\alpha}(k) + 2j - \sum_{k \geq 1} m_{\alpha}(k) 
 =  n.
\end{eqnarray*}

Now to prove Theorem \ref{th:one},
suppose $j\le n/3$. Then $n-2j\ge j-i$ for all $i\le j$, whence
$p(j-i,n-2j)=p(j-i)$ for every term in the sum that appears in
(\ref{eq:fform}), and (\ref{eq:form}) follows. If $j>n/3$ then at least the
single term with $i=0$, namely the term $p(j,n-2j)$, in the sum
(\ref{eq:fform}) is strictly less than $p(j)$, so $f(n,j)$ is strictly less
than $\sum_{i\le j}p(i)p(j-i)$, as required. $\Box$

\section{A Generating Function Equation}
We prove an identity to be used in computing $E(X_{m,i})$.
For $1 \leq i \leq m$, let
$F_{m,i}(n,k)$ be the number of partitions $\lambda$ of $n$
with $\X_{m,i}(\lambda)=k$, where $X_{m,i}$ is defined in (\ref{eq:Xdef}).
To get the generating function, note that
a part $am+b$ in the conjugate $\lambda'$,
with $1 \leq b \leq m$, contributes $a$ to the sum if $b<i$ and
$a+1$ if $b \geq i$.  So,
\begin{equation}
\sum_{n,k\geq 0} F_{m,i}(n,k) q^n u^k=
\prod_{a \geq 0}
\left( \prod_{b=1}^{i-1}\frac{1}{1-u^aq^{am+b}}
\prod_{b=i}^{m}\frac{1}{1-u^{a+1}q^{am+b}}\right).
\label{eq:gfF}
\end{equation}
Then, 
\[
\sum_{\lambda: |\lambda|=n}X_{m,i}(\lambda)=
\sum_{k\geq 0} k F_{m,i}(n,k),
\]
which can be obtained from
(\ref{eq:gfF}) 
by logarithmic
differentiation as follows.
\begin{eqnarray*}
\frac{\partial}{\partial u}\sum_{n,k \geq 0}F_{m,i}(n,k)q^nu^k \bigg|_{u=1}&=&\sum_{n\ge 0}\left(\sum_{k\geq 0} k F_{m,i}(n,k) \right)q^n\\
&=&  \prod_{j \geq 0}\frac{1}{(1-q^j)}\sum_{a \geq 0}\left( \frac{\partial}{\partial u}\left(\sum_{b=1}^{i-1}\frac{1}{1-u^aq^{am+b}}+
\sum_{b=i}^{m}\frac{1}{1-u^{a+1}q^{am+b}}\right) \right)_{u=1}\\
&=&{\cal P}(q)\sum_{a \geq 0}\left(\sum_{b=1}^{i-1}\sum_{j\geq 1}aq^{(am+b)j}+\sum_{b=i}^{m}\sum_{j\geq 1}(a+1)q^{(am+b)j}\right)\\
& = & {\cal P}(q)\sum_{j \geq 1}\sum_{d\geq 1}\left\lfloor\frac{d+m-i}{m}\right\rfloor q^{dj},
\end{eqnarray*}
in which
\[{\cal P}(z)=\prod_{j\ge 1}\frac{1}{1-z^j}\]
is Euler's unrestricted partition generating function.
It follows that
\begin{equation}
\label{eq:GF1}
\sum_{\lambda: |\lambda|=n}X_{m,i}(\lambda)=
\sum_{k\geq 1}p(n-k)\sum_{d|k} \left\lfloor\frac{d+m-i}{m}\right\rfloor.
\end{equation}
When $m=i=1$, this gives the familiar
\begin{equation}
\label{eq:GF2}
np(n)= \sum_{k\geq 1}p(n-k)\sum_{d|k}d.
\end{equation}

\section{Proof of Theorem \ref{th:three}}
\label{se:PTthree}
We start with three Lemmas.

\smallskip

\begin{lemma}
\label{le:one}
Let $m\ge 1$ and $h$ be integers.  Define
$$
\zeta_{m,h}(s)=\sum_{j=1,\,j\equiv h\, (\mod m)}^{\infty}j^{-s}, ~~~~~ \Re(s) > 1.
$$
Then,
$$
\zeta_{m,h}(s)={1/m \over s-1} + \frac{\gamma}{m} + \gamma_{m,h} + A(s),
$$
where $A(s)$ is analytic for $\Re(s)>0$; and, uniformly
for $s$ in a compact
subset of that region, we have $A(s) = O(s-1)$. 
The constants $\gamma_{m,h}$ are given by
$$
\gamma_{m,h}
= {1 \over m}\sum_{\ell=1}^{m-1}\omega^{-h\ell}\log\left({1 \over 1-\omega^{\ell}}\right),
$$
where $\omega=e^{2\pi\sqrt{-1}/m}$.
\end{lemma}

\smallskip

\noindent{\bf Proof.} There is\footnote{See also subsection \ref{subsec:one} below.} a standard technique,
sectioning, for extracting particular coefficients
from a series.  Define
$$
\hat{\zeta}_{m,\ell}(s) = \sum_{j=1}^{\infty} {\omega^{\ell j} \over j^s}.
$$
The desired result follows from the three equations
\begin{eqnarray*}
\zeta_{m,h}(s) &=& \sum_{\ell=0}^{m-1} \frac{\omega^{-\ell h}}{ m}\,\hat{\zeta}_{m,\ell}(s),\\
\hat{\zeta}_{m,0}(s) &=& \zeta(s) = {1 \over s-1} + \gamma + \O(s-1),\\
\hat{\zeta}_{m,\ell}(s) &=& \log{1 \over 1-\omega^{\ell}} + \O(s-1), ~~~
\ell \not\equiv 0\,(\mathrm{mod}\, m),
\end{eqnarray*}
where in the latter two the big-oh terms represent
functions analytic in $\Re(s)>0$.  For the displayed equation
concerning $\zeta(s)$ near $s=1$, see for example \cite{Titch}. $\Box$

\medskip

\begin{lemma}
\label{le:two}
Let $\tau(k)$ denote the divisor counting function
$$
\tau(k) = 
\sum_{d|k}1.
$$
Then, for $\Re{\alpha}>0$
$$
\sum_{k=1}^{\infty}\tau(k)e^{-k\alpha} = \alpha^{-1}\log \alpha^{-1}
+ \gamma\alpha^{-1} + \O(1).
$$
Additionally, for $t=1,2$,
$$
\sum_{k=1}^{\infty}k^t\tau(k)e^{-k\alpha} =\O\bigl( \alpha^{-t-1}
                                              \log \alpha^{-1} \bigr).
$$
\end{lemma}

\smallskip

\noindent{\bf Proof.}  Recall the Lambert series
$$
\sum_{k=1}^{\infty}\tau(k)x^k = \sum_{d=1}^{\infty}{x^d \over 1-x^d},
$$
and Mellin's formula
$$
e^{-\alpha}={1 \over 2\pi i}\int_{(a)} \alpha^{-s} \Gamma(s) ds,
~ \alpha,a>0.
$$
In the latter, $(a)$ indicates integration along the vertical
line $\Re(s)=a$.  Replace $\alpha$ by $j\alpha$, sum for
$j\ge 1$, and take $a>1$ so that the infinite sum converges:
$$
{e^{-\alpha}\over 1-e^{-\alpha}} = 
{1 \over 2\pi i}\int_{(a)} \alpha^{-s} \zeta(s) \Gamma(s) ds,
~ \alpha>0, a>1.
$$
Now, again, replace $\alpha$ by $d\alpha$ and sum for $d\ge 1$:
\begin{equation}
\label{eq:L21}
\sum_{k=1}^{\infty}\tau(k)e^{-k\alpha} = \sum_{d=1}^{\infty}{e^{-d\alpha} \over 1-e^{-d\alpha}}
 = 
{1 \over 2\pi i}\int_{(a)} \alpha^{-s} \zeta(s)^2 \Gamma(s) ds,
~ a>1.
\end{equation}
If we move the line of integration to the left, across the singularity
of the integrand at $s=1$, and use
\begin{equation}
\label{eq:L22}
{\rm Res}_{s=1}\alpha^{-s}\zeta(s)^2\Gamma(s) = \alpha^{-1}\log\alpha^{-1}
+\gamma\alpha^{-1},
\end{equation}
we obtain
$$
\sum_{k=1}^{\infty}\tau(k)e^{-k\alpha}
=  \alpha^{-1}\log\alpha^{-1}
+\gamma\alpha^{-1} +
{1 \over 2\pi i}\int_{(a)} \alpha^{-s} \zeta(s)^2 \Gamma(s) ds,
~ 0<a<1.
$$
The Lemma follows by moving the line of integration again to the
left, across the singularity of $\Gamma(s)$ at $s=0$, this time to the region $-1<a<0$.  The missing details of
this proof (for instance, interchanging the order of infinite
summation and integration, bounding integrands to justify the
movement of paths of integration) can be readily supplied by using
this fact:
along vertical lines $|\Gamma(s)|$
decays exponentially, while $|\zeta(s)|$ grows at most
polynomially \cite{Rade}.

\begin{lemma}
\label{le:three}
Let $m\ge 1$ and $h$ be integers,
$\tau_{m,h}(k)$ be the generalized divisor counting function
$$
\tau_{m,h}(k) = 
\sum_{d|k,\,d\equiv h\, (\mod m)}1,
$$
and $\gamma_{m,h}$ be the constants introduced in Lemma \ref{le:one}.
Then, for $\Re{\alpha}>0$
$$
\sum_{k=1}^{\infty}\tau_{m,h}(k)e^{-k\alpha} =
{1 \over m} \alpha^{-1}\log\alpha^{-1}
+ {1 \over m}\gamma\alpha^{-1} + \gamma_{m,h}\alpha^{-1}
+ \O(1).
$$
Additionally, for $t=1,2$,
$$
\sum_{k=1}^{\infty}k^t\tau_{m,h}(k)e^{-k\alpha} =\O\bigl( \alpha^{-t-1}\log\alpha^{-1}
                                                \bigr).
$$
\end{lemma}

\smallskip

\noindent{\bf Proof.}  The proof proceeds exactly
as the proof of Lemma \ref{le:two}.  Instead of equations (\ref{eq:L21})
and (\ref{eq:L22}) we have
$$
\sum_{k=1}^{\infty}\tau_{m,h}(k)e^{-k\alpha}
 = \sum_{d=1,d\equiv h (\mod m)}^{\infty}{e^{-d\alpha} \over 1-e^{-d\alpha}}
 = 
{1 \over 2\pi i}\int_{(a)} \alpha^{-s} \zeta_{m,h}(s)\zeta(s) \Gamma(s) ds,
~ a>1,
$$
and
$$
{\rm Res}_{s=1}\alpha^{-s}\zeta_{m,h}(s)\zeta(s)\Gamma(s) =
{1 \over m} \alpha^{-1}\log\alpha^{-1}
+ {1 \over m}\gamma\alpha^{-1} + \gamma_{m,h}\alpha^{-1},
$$
the second being obtained by an application of Lemma \ref{le:one}. 

\vskip 20pt

\noindent{\bf Proof of Theorem \ref{th:three}.}
The desired average is the sum of $\X_{m,i}(\lambda)$
over all partitions $\lambda$, divided by $\p(n)$.
{F}rom (\ref{eq:GF1})
$$
\sum_{\lambda}\X_{m,i}(\lambda) = \sum_{k=1}^n \Bigl(
                     \sum_{d|k} \left\lfloor {d+m-i \over m} \right\rfloor
                                         \Bigr) \p(n-k).
$$
We introduce the functions $\chi_j(d)$ defined by
$$
\chi_j(d) = \cases{
                      1 &  if $d \equiv j\,(\mod m)$  \cr
                      0 &  otherwise,
                  }
$$
and note
$$
\left\lfloor {d+m-i \over m} \right\rfloor = {d+m-i \over m} -
\sum_{j=1}^{m-1} {j \chi_{i+j}(d) \over m}.
$$
Since by (\ref{eq:GF2})
$$
n\p(n) = \sum_{k=1}^n \Bigl( \sum_{d|k}d \Bigr) \p(n-k),
$$
we have
\begin{equation}
\label{eq:SumFormula}
\sum_{\lambda}\X_{m,i}(\lambda) = {n \over m}\p(n) +
{m-i \over m}S - \sum_{j=1}^{m-1} {j \over m} S_{i+j},
\end{equation}
where
$$
S = \sum_{k=1}^n \tau(k) \p(n-k)
$$
and
$$
S_j = \sum_{k=1}^n \tau_{m,j}(k) \p(n-k).
$$
The plan of the proof is to obtain asymptotic estimations of
$S$ and $S_j$, and then substitute these into equation (\ref{eq:SumFormula}).

To estimate $S$, we replace $\p(n-k)$ from
$$
\p(n-k) = \p(n) e^{-Ck/(2\sqrt{n})} \Bigl(
           1 + \O(k/n+k^2/n^{3/2})
                                        \Bigr),
$$
uniformly for $k\le n^{2/3}$.
(This is an easy consequence of the Hardy-Ramanujan formula for
$\p(n)$.)  We have
\begin{eqnarray*}
S &=& \sum_{1 \le k \le n^{2/3}} \tau(k)\p(n-k) +
\sum_{n^{2/3} < k \le n} \tau(k)\p(n-k) \\
  &=& \p(n) \sum_{1 \le k \le n^{2/3}} \tau(k)
e^{-Ck/(2\sqrt{n})} \Bigl(
           1 + \O(k/n+k^2/n^{3/2})
                                        \Bigr)
+ \O\bigl(\p(n)n^2e^{-(C/2)n^{1/6}}
    \bigr) \\
&=& \p(n) \sum_{k=1}^{\infty} \tau(k)
e^{-Ck/(2\sqrt{n})} \Bigl(
           1 + \O(k/n+k^2/n^{3/2})
                                        \Bigr)
+ \O\bigl(\p(n)n^2e^{-(C/2)n^{1/6}}
    \bigr),
\end{eqnarray*}
where we have used the extremely crude bound
$\tau(k)\le k$ at two points.
The infinite sum appearing here can be estimated
by Lemma \ref{le:two}, with $\alpha=C/(2\sqrt{n})$.  We
conclude that
$$
S = \p(n) \left(
\frac{1}{C}\sqrt{n} \log n + \left(\frac{2}{C}(\gamma+\log{\frac{2}{C}})\right)\sqrt{n}
+ \O(\log n)
          \right).
$$ 
The estimation of $S_j$ proceeds very similarly, using
Lemma \ref{le:three} instead of Lemma \ref{le:two}.  The
bottom line is
\[
S_j =\frac{p(n)}{mC} \left(
\sqrt{n} \log n +\left(2\log{\left(\frac{2}{C}\right)}+2\gamma+2m\gamma_{m,j}\right)\sqrt{n}+ \O(\log n)
          \right).
\]
When
we substitute the formulas for $S$ and $S_j$ into equation
(\ref{eq:SumFormula}), and perform some elementary algebra
including
$$
{m-i \over m} - {1 \over m}\sum_{j=1}^{m-1}{j \over m} = {m+1-2i \over 2m},
$$
and
$$
\sum_{j=1}^{m-1} {j \over m}\gamma_{m,i+j} = -{1 \over m}\sum_{\ell=1}^{m-1}
{\omega^{-\ell(i-1)} \over 1-\omega^\ell}
   \log(1-\omega^{\ell}),
$$
(where as usual $\omega=e^{2\pi\sqrt{-1}/m}$),
we obtain the theorem.

\section{Notes and Bibliographic Remarks}
\subsection{Remarks on Lemma \ref{le:one}}
\label{subsec:one}
Our $\zeta_{m,h}(s)$ is related to both the L-functions of Dirichlet and the generalized zeta function of Hurwitz, defined by 
\[\zeta(s,a)=\sum_{n\ge 0}\frac{1}{(n+a)^s}.\]
Indeed, our $\zeta_{m,h}(s)=m^{-s}\zeta(s,h/m)$. It is known (e.g., \cite{ww}, Chap. XIII) that
\[\zeta(s,a)=\frac{1}{s-1}-\frac{\Gamma'(a)}{\Gamma(a)}+O(s-1)\qquad (s\to 1),\]
in which the ``big-oh'' term is an entire function. Consequently, except possibly for the given form of the constants $\gamma_{m,h}$, Lemma \ref{le:one} is well known.
\subsection{Remarks on Lemma \ref{le:two}}
The asymptotics of the Lambert series whose coefficients are the divisor function are due to Wigert \cite{wi}, and are used in \cite{Titch} to investigate the mean values of the Riemann Zeta function on the critical line. We included the derivation here because it applies equally well to Lemma \ref{le:three}, which is more general. In fact, the process of moving the line of integration to the left by one unit can be indefinitely repeated and the result is the complete asymptotic series representation for this Lambert series (see, e.g.,  \cite{be}). Although this asymptotic series is well known for the standard divisor function, we have not been able to find in the literature the corresponding result for the ``mod $m$'' divisor function $\tau_{m,h}(k)$. Since this is an easy consequence of our method, we quote the final result here, which is the following asymptotic series:
\begin{equation}
\sum_{k=1}^{\infty}\tau_{m,h}(k)e^{-k\alpha}\approx \frac{1}{m}\alpha^{-1}\log{\alpha^{-1}}+\left(\frac{\gamma}{m}+\gamma_{m,h}\right)\alpha^{-1}-\sum_{n=0}^{\infty}\frac{B_{n+1}B_{n+1}\kern-4pt\left(\frac{h}{m}\right)}{(n+1)!(n+1)}(\alpha m)^n,
\end{equation}
in which both $B_n$, the $n$th Bernoulli \textit{number}, and $B_{n}(x)$, the $n$th Bernoulli \textit{polynomial}, appear. This result agrees with known (\cite{be}, \cite{wi}) formulas for the usual divisor function in the case $m=h=1$.

This use of Mellin's formula is found frequently in analytic number theory.  Another example is seen
in Chapter 6 of \cite{Andrews}, which deals with the asymptotics of coefficients of infinite products.  The method was also used in Husimi's \cite{Husimi} 1938 paper on the average number of parts in a random partition.
\subsection{Remarks on the constants $\gamma_{m,h}$}
In Lemma \ref{le:one} we have given an explicit formula for the constants $\gamma_{m,h}$. Here we note that a purely real form for these constants had previously been found by Gauss, namely
\begin{equation}
\label{eq:gssform}
\gamma_{m,h}
= \frac{1}{m}\left(\frac{\pi}{2}\cot{\left(\frac{h}{m}\pi\right)}+\log{2}
-2\kern-5pt\sum_{0<k<m/2}\cos{\left(\frac{2hk}{m}\pi\right)}\cdot\log{\sin{\frac{k\pi}{m}}}\right).
\end{equation}

Indeed, as we mentioned earlier, our $\zeta_{m,h}(s)=m^{-s}\zeta(s,h/m)$, where the $\zeta$ on the right is the Hurwitz Zeta function. It is known (e.g., \cite{ww}, Chap. XIII) that
\[\zeta(s,a)=\frac{1}{s-1}-\frac{\Gamma'(a)}{\Gamma(a)}+O(s-1)\qquad (s\to 1),\]
in which the ``big-oh'' term is an entire function.
Hence,
\begin{eqnarray*}
\zeta_{m,h}(s)&=&m^{-s}\left(\frac{1}{s-1}-\frac{\Gamma'(\frac{h}{m})}{\Gamma(\frac{h}{m})}+O(s-1)\right)\\
&=&\frac{1}{m}\frac{1}{s-1}-\frac{\log{m}}{m}-\frac{\Gamma'(\frac{h}{m})}{m\Gamma(\frac{h}{m})}+O(s-1).
\end{eqnarray*}
Now the logarithmic derivative of the $\Gamma$-function evaluated at a rational argument can be expressed in a nice form. We have first, from section 12.16 of \cite{ww}, the relation 
\[\frac{d}{dz}\log{\Gamma(z)}=-\gamma-\frac{1}{z}+z\sum_{n\ge 1}\frac{1}{n(n+z)},\]
and therefore
\begin{eqnarray*}
 \frac{\Gamma'(\frac{h}{m})}{\Gamma(\frac{h}{m})}&=&-\gamma-\frac{m}{h}+\frac{h}{m}\sum_{n\ge 1}\frac{1}{n(n+\frac{h}{m})}\\
&=&-\gamma-\sum_{n\ge 0}\left(\frac{1}{n+\frac{h}{m}}-\frac{1}{n+1}\right).
\end{eqnarray*}

Next, from Ex. 1.2.9.18 of \cite{kn}, we have the following identity, which is attributed to Gauss:
\begin{equation}
\label{eq:gauss}
\sum_{n\ge 0}\left(\frac{1}{n+p/q}-\frac{1}{n+1}\right)=\frac{\pi}{2}\cot{\frac{p}{q}\pi}+\log{2q}-2\sum_{0<k<q/2}\cos{\frac{2pk}{q}\pi}\cdot\log{\sin{\frac{k\pi}{q}}}.
\end{equation}
Thus we have\footnote{This is also due to Gauss. See \cite{al}.}, for $0<h<m$,
\begin{equation}
\frac{\Gamma'(\frac{h}{m})}{\Gamma(\frac{h}{m})}=-\gamma -\frac{\pi}{2}\cot{\left(\frac{h}{m}\pi\right)}-\log{(2m)}+2\sum_{0<k<m/2}\cos{\left(\frac{2hk}{m}\pi\right)}\cdot\log{\sin{\frac{k\pi}{m}}}.
\end{equation}
This shows that $\zeta_{m,h}(s)$ is of the form
\[\frac{1}{m}\frac{1}{s-1}+\frac{1}{m}\left(\gamma +\frac{\pi}{2}\cot{\left(\frac{h}{m}\pi\right)}+\log{2}
-2\sum_{0<k<m/2}\cos{\left(\frac{2hk}{m}\pi\right)}\cdot\log{\sin{\frac{k\pi}{m}}}\right)+O(s-1),\]
which is the result (\ref{eq:gssform}). $\Box$

\end{document}